\documentclass[11pt,leqno]{article}
\usepackage{amsmath, amsthm, amssymb}
\usepackage{amsfonts}
\usepackage[algo2e,boxed]{algorithm2e}

\newcommand{\ignore}[1]{}
\ignore{The text in between these parentheses is ignored by LaTeX
ignore}

\def\@begintheorem#1#2{\par\bgroup{\sc #1\ #2. }\it\ignorespaces}
\def\@opargbegintheorem#1#2#3{\par\bgroup{\sc #1\ #2\ (#3). } \it\ignorespaces}
\def\@endtheorem{\egroup}
\newtheorem{theorem}{Theorem}[section]
\newtheorem{corollary}[theorem]{Corollary}
\newtheorem{lemma}[theorem]{Lemma}
\newtheorem{proposition}[theorem]{Proposition}
\newtheorem{problem}[theorem]{Problem}
\newtheorem{example}[theorem]{Example}
\newtheorem{algorithm}[theorem]{Algorithm}
\newtheorem{definition}[theorem]{Definition}
\newcommand{\bt}[1]{\begin{theorem}\label{#1}}
\newcommand{\bc}[1]{\begin{corollary}\label{#1}}
\newcommand{\bl}[1]{\begin{lemma}\label{#1}}
\newcommand{\bp}[1]{\begin{proposition}\label{#1}}
\newcommand{\bpro}[1]{\begin{problem}\label{#1}}
\newcommand{\be}[1]{\begin{example}\rm\label{#1}}
\newcommand{\ba}[1]{\begin{algorithm}\rm\label{#1}}
\newcommand{\bd}[1]{\begin{definition}\rm\label{#1}}
\newcommand{\bpr}{\begin{proof}}
\newcommand{\et}{\end{theorem}}
\newcommand{\ec}{\end{corollary}}
\newcommand{\el}{\end{lemma}}
\newcommand{\ep}{\end{proposition}}
\newcommand{\epro}{\end{problem}}
\newcommand{\ee}{\end{example}}
\newcommand{\ea}{\end{algorithm}}
\newcommand{\ed}{\end{definition}}
\newcommand{\epr}{\end{proof}}

\def\R{\mathbb{R}}
\def\Z{\mathbb{Z}}
\def \supp {{\rm supp}}
\def \fr {{\rm F}}
\def \la {\lambda}

\begin{document}

\title{\bf Nonlinear Optimization over a\\ Weighted Independence System}
\author{Jon Lee \and Shmuel Onn \and Robert Weismantel}
\date{}
\maketitle

\begin{abstract}
We consider the problem of optimizing a nonlinear objective
function over a weighted independence system presented by a
linear-optimization oracle. We provide a polynomial-time algorithm
that determines an r-best solution for nonlinear functions of the
total weight of an independent set, where r is a constant that
depends on certain Frobenius numbers of the individual weights
and is independent of the size of the ground set. In contrast,
we show that finding an optimal (0-best) solution requires
exponential time even in a very special case of the problem.
\end{abstract}

\section{Introduction}

An \emph{independence system} is a nonempty set  of vectors
$S\subseteq\{0,1\}^n$ with the property that $x\in\{0,1\}^n$~,
$x\leq y\in S$ implies $x\in S$~. The general nonlinear optimization
problem over a multiply-weighted independence system is as follows.

\vskip.2cm\noindent {\bf Nonlinear optimization over a
multiply-weighted independence system.} Given independence system
$S\subseteq\{0,1\}^n$~, weight vectors $w^1,\dots,w^d\in\Z^n$~, and
function $f:\Z^d\rightarrow\R$~, find $x\in S$ minimizing the
objective
\[
f(w^1x,\dots,w^dx)\ =\ f\left(\sum_{j=1}^n
w^1_jx_j,\dots,\sum_{j=1}^n w^d_jx_j\right)~.
\]

\vskip.2cm\noindent The representation of the objective in the above
composite form has several advantages. First, for $d>1$~, it can
naturally be interpreted as \emph{multi-criteria optimization}: the
$d$ given weight vectors $w^1,\dots,w^d$ represent $d$ different
criteria, where the value of $x\in S$ under criterion $i$ is its
$i$-th total weight $w^ix=\sum_{j=1}^n w^i_jx_j$~; and the objective
is to minimize the ``balancing" $f(w^1x,\dots,w^dx)$ of the $d$
given criteria by the given function $f$~. Second, it allows us to
classify nonlinear optimization problems into a hierarchy of
increasing generality and complexity: at the bottom lies standard
linear optimization, recovered with $d=1$ and $f$ the identity on
$\Z$~; and at the top lies the problem of minimizing an arbitrary
function, which is typically intractable, arising with $d=n$ and
$w_i={\bf 1}_i$ the $i$-th standard unit vector in $\Z^n$ for all
$i$~.

The computational complexity of the problem depends on the number
$d$ of weight vectors, on the weights $w^i_j$~, on the type of
function $f$ and its presentation, and on the type of independence
system $S$ and its presentation. For example, when $S$ is a
\emph{matroid}, the problem can be solved in polynomial time for any
fixed $d$~, any $\{0,1,\dots,p\}$-valued weights $w^i_j$ with $p$
fixed, and any function $f$ presented by a \emph{comparison oracle},
even when $S$ is presented by a mere \emph{membership oracle}, see
\cite{BLMORWW}. Also, when $S$ consists of the \emph{matchings} in a
given bipartite graph $G$~, the problem can be solved in polynomial
time for any fixed $d$~, any weights $w^i_j$ presented in unary, and
any \emph{concave} function $f$~, see \cite{BO}; but on the other
hand, for \emph{convex} $f$~, already with fixed $d=2$ and
$\{0,1\}$-valued weights $w^i_j$~, it includes as a special case the
notorious \emph{exact matching problem}, the complexity of which is
long open \cite{MVV,PY}.

In view of the difficulty of the problem already for $d=2$~, in this
article we take a first step and concentrate on \emph{nonlinear
optimization over a (singly) weighted independence system}, that is,
with $d=1$~, single weight vector $w=(w_1,\dots,w_n)\in\Z^n$~, and
univariate function $f:\Z\rightarrow\R$~. The function $f$ can be
arbitrary and is presented by a \emph{comparison oracle} that,
queried on $x,y\in\Z$~, asserts whether or not $f(x)\leq f(y)$~. The
weights $w_j$ take on values in a $p$-tuple $a=(a_1,\dots,a_p)$ of
positive integers. Without loss of generality we assume that
$a=(a_1,\dots,a_p)$ is \emph{primitive}, by which we mean that the
$a_i$ are distinct positive integers having greatest common divisor
$\gcd(a):=\gcd(a_1,\dots,a_p)$ that is equal to $1$~. The independence
system $S$ is presented by a \emph{linear-optimization oracle} that,
queried on vector $v\in\Z^n$~, returns an element $x\in S$ that
maximizes the linear function $vx=\sum_{j=1}^n v_jx_j$~.
It turns out that solving this problem to optimality
may require exponential time (see Theorem \ref{lower_bound}),
and so we settle for an approximate solution in the
following sense, that is interesting in its own right. For a
nonnegative integer $r$~, we say that $x^*\in S$ is an
\emph{$r$-best solution} to the optimization problem over $S$ if
there are at most $r$ better objective values attained by feasible
solutions. In particular, a $0$-best solution is optimal. Recall
that the \emph{Frobenius number} of a primitive $a$ is the largest
integer $\fr(a)$ that is not expressible as a nonnegative integer
combination of the $a_i$~. We prove the following theorem.

\bt{Main} For every primitive $p$-tuple $a=(a_1,\dots,a_p)$~, there
is a constant $r(a)$ and an algorithm that, given any independence
system $S\subseteq\{0,1\}^n$ presented by a linear-optimization
oracle, weight vector $w\in\{a_1,\dots,a_p\}^n$~, and function
$f:\Z\rightarrow\R$ presented by a comparison oracle, provides an
$r(a)$-best solution to the nonlinear problem $\min\{f(wx)~:~x\in
S\}$~, in time polynomial in $n$~. Moreover:
\begin{enumerate}
\item
If $a_i$ divides $a_{i+1}$ for $i=1,\dots, p-1$~, then the algorithm
provides an optimal solution.
\item
For $p=2$~, that is, for $a=(a_1,a_2)$~, the algorithm provide an
$\fr(a)$-best solution.
\end{enumerate}
\et
In fact, we give an explicit upper bound on $r(a)$ in terms of
the Frobenius numbers of certain subtuples derived from $a$~.

Because $F(2,3)=1$~, Theorem \ref{Main} (Part 2) assures us that we can
efficiently compute a $1$-best solution in that case. It is natural
to wonder then whether, in this case, an optimal (i.e., $0$-best) solution
can be calculated in polynomial time. The next result indicates that
this cannot be done.

\bt{Main2}
There is no polynomial time algorithm for computing
an optimal (i.e., $0$-best) solution of the nonlinear optimization
problem $\min\{f(wx)\,:\,x\in S\}$ over an independence system
presented by a linear optimization oracle with $f$ presented by a
comparison oracle and weight vector $w\in\{2,3\}^n$.
\et

The next sections gradually develop the various necessary
ingredients used to establish our main results. \S\ref{Notation} sets some
notation. \S\ref{Naive} discusses a na\"{\i}ve solution strategy
that does not directly lead to a good approximation, but is a basic
building block that is refined and repeatedly used later on.
\S\ref{Partitioning} describes a way of partitioning an independence
system into suitable pieces, on each of which a suitable refinement
of the na\"{\i}ve strategy will be applied separately.
\S\ref{Monoids} provides some properties of monoids and Frobenius
numbers that will allows us to show that the refined na\"{\i}ve
strategy applied to each piece gives a good approximation within
that piece. \S\ref{Algorithm} combines all ingredients developed in
\S\ref{Naive}--\ref{Monoids}, provides a bound on the approximation
quality $r(a)$~, and provides the algorithm establishing Theorem
\ref{Main}. \S\ref{LowerBound} demonstrates that finding an optimal
solution is provably intractable, proving a refined version of
Theorem \ref{Main2}. \S\ref{Conclusion} concludes with
some final remarks and questions.

\section{Some Notation}\label{Notation}

In this section we provide some notation that will be used
throughout the article. Some more specific notation will be
introduced in later sections. We denote by $\R$~, $\R_+$~, $\Z$ and
$\Z_+$~, the reals, nonnegative reals, integers and nonnegative
integers, respectively. For a positive integer $n$, we let
$N:=\{1,\dots,n\}$~. The $j$-th standard unit vector in $\R^n$ is
denoted by ${\bf 1}_j$~. The \emph{support} of $x\in\R^n$ is the
index set $\supp(x):=\{j~:~x_j\neq 0\}\subseteq N$ of nonzero
entries of $x$~. The \emph{indicator} of a subset $J\subseteq N$ is
the vector ${\bf 1}_J:=\sum_{j\in J}{\bf 1}_j\in\{0,1\}^n$~, so that
$\supp({\bf 1}_J)=J$~. The \emph{positive} and \emph{negative} parts
of a vector $x\in\R^n$ are denoted, respectively, by
$x^+,x^-\in\R_+^n$~, and defined by $x^+_i:=\max\{x_i,0\}$ and
$x^-_i:=-\min\{x_i,0\}$ for $i=1,\dots,n$~. So, $x=x^+-x^-$~, and
$x^+_i x^-_i = 0$ for $i=1,\ldots,n$~.

Unless otherwise specified, $x$ denotes an element of $\{0,1\}^n$
and $\la,\mu,\tau,\nu$ denote elements of $\Z^p_+$~. Throughout,
$a=(a_1,\dots,a_p)$ is a \emph{primitive} $p$-tuple, by which we
mean that the $a_i$ are distinct positive integers having greatest
common divisor $\gcd(a):=\gcd(a_1,\dots,a_p)$ equal to $1$~. We will
be working with weights taking values in $a$~, that is, vectors
$w\in\{a_1,\dots,a_p\}^n$~. With such a weight vector $w$ being
clear from the context, we let $N_i:=\{j\in N~:~w_j=a_i\}$ for
$i=1,\dots,p$~, so that $N=\biguplus_{i=1}^p N_i$~. For
$x\in\{0,1\}^n$ we let $\la_i(x):=|\supp(x)\cap N_i|$ for
$i=1,\dots,p$~, and $\la(x):=(\la_1(x),\dots,\la_p(x))$~, so that
$wx=\la(x)a$~. For integers $z,s\in\Z$ and a set of integers
$Z\subseteq \Z$~, we define $z+sZ:=\{z+sx~:~x\in Z\}$~.

\section{A Na\"{\i}ve Strategy}\label{Naive}

Consider a set $S\subseteq\{0,1\}^n$~, weight vector
$w\in\{a_1,\dots,a_p\}^n$~, and function $f:\Z\rightarrow\R$
presented by a comparison oracle. Define the \emph{image} of $S$
under $w$ to be the set of values $wx$ taken by elements of $S$~,
\[\textstyle
w\cdot S\quad :=\quad \left\{wx=\sum_{j=1}^nw_jx_j ~:~ x\in
S\right\}\quad\subseteq\quad\Z_+~.
\]
As explained in the introduction, for a nonnegative integer $r$~, we
say that $x^*\in S$ is an $r$-best solution if there are at most $r$
better objective values attained by feasible solutions. Formally,
$x^*\in S$ is an \emph{$r$-best solution} if
\[
\left| \left\{ f(wx) ~:~ f(wx)< f(wx^*)~,~ x\in S \right\} \right|
~\le~ r\ .
\]
We point out the following simple observation.

\bp{image} If $f$ is given by a comparison oracle, then a necessary
condition for any algorithm to find an $r$-best solution to the
problem $\min\{f(wx)~:~x\in S\}$ is that it computes all but at most
$r$ values of the image $w\cdot S$ of $S$ under $w$~. \ep

Note that this necessary condition is also sufficient for computing
the weight $wx^*$ of an $r$-best solution, but not for computing an
actual $r$-best solution $x^*\in S$~, which may be harder.

Any point ${\bar x}$ attaining $\max\{wx~:~x\in S\}$ provides an
approximation of the image given by
\begin{equation}\label{approximation}
\{wx ~:~ x\leq {\bar x}\} \quad \subseteq \quad w\cdot S \quad
\subseteq \quad \{0,1,\dots,w{\bar x}\}\ .
\end{equation}
This suggests the following natural na\"{\i}ve strategy for finding
an approximate solution to the optimization problem over an
independence system $S$ that is presented by a linear-optimization
oracle.

\vskip.2cm
\begin{algorithm2e}[H]
\vskip.2cm \label{naive-strategy}

\begin{center}
{\bf Na\"{\i}ve Strategy}
\end{center}

{\bf input} independence system $S\subseteq\{0,1\}^n$ presented by a
linear-optimization oracle, $f:\Z\rightarrow\R$ presented by a
comparison oracle, and  $w\in\{a_1,\ldots,a_p\}^n$~\;

{\bf obtain }${\bar x}$ attaining $\max\{wx~:~x\in S\}$ using the
linear-optimization oracle for $S$~\;

{\bf  output} $x^*$ as one attaining $\min\{f(wx)~:~x\leq {\bar x}\}$
using the algorithm of Lemma \ref{strategy-lemma} below\ .

\end{algorithm2e}

\noindent Unfortunately, as the next example shows, the number of
values of the image that are missing from the approximating set on
the left-hand side of equation (\ref{approximation}) cannot
generally be bounded by any constant. So by Proposition \ref{image},
this strategy cannot be used \emph{as is} to obtain a provably good
approximation.

\be{example} Let $a:=(1,2)$~, $n:=4m$~, $y:=\sum_{i=1}^{2m}{\bf
1}_i$~, $z:=\sum_{i=2m+1}^{4m}{\bf 1}_i$~, and $w:=y+2z$~, that is,
$$y \ =\ (1,\dots,1,0,\dots,0)\,,\quad z \ =\ (0,\dots,0,1,\dots,1)\,,\quad
w \ =\ (1,\dots,1,2,\dots,2)\,,$$ define $f$ on $\Z$ by
\[
f(k):= \left\{
  \begin{array}{ll}
    k~, & \hbox{$k$ odd;} \\
    2m~, & \hbox{$k$ even,}
  \end{array}
\right.
\]
and let $S$ be the independence system
\[
S\ :=\ \{x\in\{0,1\}^n ~:~ x\leq y\}\ \cup\ \{x\in\{0,1\}^n ~:~
x\leq z\}~.
\]
Then the unique optimal solution of the linear-objective problem
$\max\{wx ~:~ x\in S\}$ is ${\bar x}:=z$~, with $w{\bar x}=4m$~, and
therefore
\begin{eqnarray*}
&&\{wx ~:~ x\leq {\bar x}\}\ =\ \{2i ~:~ i=0,1,\dots,2m\}~, \mbox{ and}\\
&&w\cdot S \ =\ \{i ~:~ i=0,1,\dots,2m\}\ \cup\ \{2i ~:~
i=0,1,\dots,2m\}~.
\end{eqnarray*}
So all $m$ odd values (i.e., $1,3,\ldots,2m-1$) in the image $w\cdot
S$ are missing from the approximating set $\{wx ~:~ x\leq {\bar
x}\}$ on the left-hand side of (\ref{approximation}), and $x^*$
attaining $\min\{f(wx) ~:~ x\leq{\bar x}\}$ output by the above
strategy has objective value $f(wx^*)=2m$~, while there are
$m={n\over 4}$ better objective values (i.e., $1,3,\ldots,2m-1$)
attainable by feasible points (e.g., $\sum_{i=1}^k {\bf 1}_i$~, for
$k=1,3,\ldots,2m-1$). \ee \vskip0.2in

Nonetheless, a more sophisticated refinement of the na\"{\i}ve
strategy, applied repeatedly to several suitably chosen subsets of
$S$ rather than $S$ itself, will lead to a good approximation. In
the next two sections, we develop the necessary ingredients that
enable us to implement such a refinement of the na\"{\i}ve strategy
and to prove a guarantee on the quality of the approximation it
provides. Before proceeding to the next section, we note that the
na\"{\i}ve strategy can be efficiently implemented as follows.

\bl{strategy-lemma} For every fixed $p$-tuple $a$~, there is a
polynomial-time algorithm that, given univariate function
$f:\Z\rightarrow\R$ presented by a comparison oracle, weight vector
$w\in\{a_1,\dots,a_p\}^n$~, and $\bar x\in\{0,1\}^n$~, solves
\begin{equation*}
\min\{f(wx) ~:~ x\leq {\bar x}\}\ .
\end{equation*}
\el \bpr Consider the following algorithm:

\begin{algorithm2e}[H]
\label{strategy-algorithm}

{\bf input} function $f:\Z\rightarrow\R$ presented by a comparison
oracle, $w\in\{a_1,\ldots,a_p\}^n$ and $\bar x\in\{0,1\}^n$~\;

{\bf let} $N_i:=\{j ~:~ w_j=a_i\}$ and $\tau_i:=\la_i({\bar
x})=|\supp({\bar x})\cap N_i|,\,$ $i=1,\dots,p$~\;

\For {\emph{every choice of}
$\nu=(\nu_1,\dots,\nu_p)\leq(\tau_1,\dots,\tau_p)=\tau$}{

{\bf  determine} some $x_\nu\leq{\bar x}$ with
$\la_i(x_\nu)=|\supp(x_\nu)\cap N_i|=\nu_i,\,$ $i=1,\dots,p$~\;}

{\bf  output} $x^*$ as one minimizing $f(wx)$ among the $x_\nu$ by
using the comparison oracle of $f$~\ .

\end{algorithm2e}

\noindent Since the value $wx$ depends only on the cardinalities
$|\supp(x)\cap N_i|,\,$ $i=1,\dots,p$~, it is clear that
\[
\{wx ~:~ x\leq {\bar x}\}\quad =\quad \{wx_\nu ~:~ \nu\leq\tau\}~.
\]
Clearly, for each choice $\nu\leq\tau$ it is easy to determine some
$x_\nu\leq{\bar x}$ by zeroing out suitable entries of $\bar x$~.
The number of choices $\nu\leq\tau$ and hence of loop iterations and
comparison-oracle queries of $f$ to determine $x^*$ is
\[
\prod_{i=1}^p (\tau_i+1)\ \leq\ (n+1)^p~.
\]
\epr

\section{Partitions of Independence Systems}
\label{Partitioning}

Define the \emph{face of $S\subseteq\{0,1\}^n$ determined by two
disjoint subsets $L,U\subseteq N=\{1,\ldots,n\}$} to be
\[
S^U_L ~:=~ \{x\in S ~:~ x_j=0\ \mbox{for}\ j\in L\,,\ x_j=1\
\mbox{for}\ j\in U\}~.
\]
Our first simple lemma reduces linear optimization over faces of $S$
to linear optimization over $S$~. \bl{face-lemma} Consider any
nonempty set $S\subseteq\{0,1\}^n$~, weight vector $w\in\Z^n$~, and
disjoint subsets $L,U\subseteq N$~. Let $\alpha:=1+2n\max|w_j|$~,
let ${\bf 1}_L,{\bf 1}_U\in\{0,1\}^n$ be the indicators of $L,U$
respectively, and let
\begin{eqnarray}\label{face_equation}
\nonumber
v &:=& \max\left\{(w+\alpha({\bf 1}_U-{\bf 1}_L))x ~:~ x\in S\right\} ~-~ |U|\alpha \\
  & = & \textstyle \max\left\{ wx - \alpha\left( \sum_{j\in U} (1-x_j) + \sum_{j\in L} x_j
     \right) ~:~ x\in S \right\}~.
\end{eqnarray}
Then either $v>-{1\over2}\alpha$, in which case $\max\{wx ~:~ x\in
S^U_L\}=v$ and the set of maximizers of $wx$ over $S^U_L$ is equal
to the set of maximizers of the program (\ref{face_equation}), or
$v<-{1\over2}\alpha$~, in which case $S^U_L$ is empty. \el \bpr For
all $x\in\{0,1\}^n$~, we have $-{1\over2}\alpha < wx <
{1\over2}\alpha$~, and so for all $y\in S\setminus S^U_L$ and $z\in
S^U_L$ we have
\begin{eqnarray*}
wy - \alpha\left( \sum_{j\in U} (1-y_j) + \sum_{j\in L} y_j  \right)
& \le & wy - \alpha ~<~ {1\over2}\alpha -\alpha ~=~ - {1\over2}\alpha \\
& <   & wz \ =\  wz - \alpha\left( \sum_{j\in U} (1-z_j) +
\sum_{j\in L} z_j  \right)~.
\end{eqnarray*}
\epr

Let $S\subseteq\{0,1\}^n$ and $w\in\{a_1,\dots,a_p\}^n$ be
arbitrary, and let $N_i:=\{j\in N ~:~ w_j = a_i\}$ as usual. As
usual, for $x\in S$~, let
 $\la_i(x):=|\supp(x)\cap N_i|$ for each $i$~.
For $p$-tuples $\mu=(\mu_1,\dots,\mu_p)$ and $\la =
(\la_1,\dots,\la_p)$ in $\Z^p_+$ with $\mu\le \la$~, define
\begin{equation}\label{partition-equation}
S^\la_\mu\ :=\ \left\{x\in S ~:~
  \begin{array}{ll}
     \la_i(x)=\mu_i~, & \hbox{if $\mu_i<\la_i$}~, \\
   \la_i(x)\geq\mu_i~, & \hbox{if $\mu_i=\la_i$~.}
  \end{array}
\right\}~.
\end{equation}

\bp{Partition} Let $S\subseteq\{0,1\}^n$ be arbitrary. Then every
$\la\in\Z^p_+$ induces a partition of $S$ given by
$$S\ =\ \biguplus_{\mu\leq\la} S^\la_\mu\ .$$
\ep \bpr Consider any $x\in S$~, and define $\mu\leq\la$ by $\mu_i
:= \min\{\la_i(x),\la_i\}$~.
Then  $x\in S^\la_\mu$~, but $x\notin S^\la_\nu$ for $\nu\leq\la$~,~
$\nu\neq\mu$~. \epr

\bl{mu-lemma} For all fixed $p$-tuples $a$ and $\la\in\Z^p_+$~,
there is a polynomial-time algorithm that, given any independence
system $S$ presented by a linear-optimization oracle,
$w\in\{a_1,\dots,a_p\}^n$~, and $\mu\in\Z^p_+$ with $\mu\leq\la$~,
solves
\begin{equation*}
\max\left\{wx ~:~ x\in S^\la_\mu\right\}~.
\end{equation*}
\el \bpr Consider the following algorithm:

\begin{algorithm2e}[H]
\label{mu-algorithm}

{\bf input} independence system $S\subseteq\{0,1\}^n$ presented by a
linear-optimization oracle~, $w\in\{a_1,\ldots,a_p\}^n$~, and
$\mu\leq\la$~\;

{\bf let} $I:=\{i ~:~ \mu_i<\la_i\}$ and $N_i:=\{j\in N ~:~
w_j=a_i\},\,$ $i=1,\dots,p$~\;

\For{ \emph{every} $S_i\subseteq N_i$ \emph{with} $|S_i|=\mu_i,\,$
$i=1,\dots,p,\,$ \emph{if any,}}{

{\bf let} $L:=\bigcup_{i\in I}\left(N_i\setminus S_i\right)$
                         and $U:=\bigcup_{i=1}^p S_i$~\;

{\bf find} by the algorithm of Lemma \ref{face-lemma} an
$x(S_1,\dots,S_p)$ attaining $\max\{wx ~:~ x\in S^U_L\}$ if any\;

}

{\bf  output} $x^*$ as one maximizing $wx$ among all of the
$x(S_1,\dots,S_p)$ (if any) found in the loop above\ .

\end{algorithm2e}

\noindent It is clear that $S^\la_\mu$ is the union of the $S^U_L$
over all choices $S_1,\dots,S_p$ as above, and therefore $x^*$ is
indeed a maximizer of $wx$ over $S^\la_\mu$~. The number of such
choices and hence of loop iterations is
\[
\prod_{i=1}^p {|N_i|\choose \mu_i}\ \leq\ \prod_{i=1}^p n^{\mu_i} \
\leq\ \prod_{i=1}^p n^{\la_i}~,
\]
which is polynomial because $\la$ is fixed. In each iteration, we
find $x(S_1,\dots,S_p)$ maximizing $wx$ over $S^U_L$ or detect
$S^U_L=\emptyset$ by applying the algorithm of Lemma
\ref{face-lemma} using a single query of the linear-optimization
oracle for $S$~. \epr

We will later show that, for a suitable choice of $\la$~, we can
guarantee that, for every block $S^\la_\mu$ of the partition of $S$
induced by $\la$~, the na\"{\i}ve strategy applied to $S^\la_\mu$
does give a good solution, with only a constant number of better
objective values obtainable by solutions within $S^\la_\mu$~. For
this, we proceed next to take a closer look at the monoid generated
by a $p$-tuple $a$ and at suitable restrictions of this monoid.

\section{Monoids and Frobenius Numbers}\label{Monoids}

Recall that a $p$-tuple $a=(a_1,\dots,a_p)$ is \emph{primitive} if
the $a_i$ are distinct positive integers having greatest common
divisor $\gcd(a)=\gcd(a_1,\dots,a_p)$ is $1$~. For $p=1$~, the only
primitive $a=(a_1)$ is the one with $a_1=1$~. The \emph{monoid} of
$a=(a_1,\dots,a_p)$ is the set of nonnegative integer combinations
of its entries,
\[\textstyle
M(a)\ =\ \left\{\mu a=\sum_{i=1}^p\mu_i a_i ~:~ \mu\in
\Z_+^p\right\}~.
\]
The \emph{gap set} of $a$ is the set $G(a):=\Z_+\setminus M(a)$ and
is well known to be finite \cite{Bra}. If all $a_i\geq 2$~, then
$G(a)$ is nonempty, and its maximum element is known as the
\emph{Frobenius number} of $a$~, and will be denoted by
$\fr(a):=\max G(a)$~. If some $a_i=1$~, then $G(a)=\emptyset$~, in
which case we define $\fr(a):=0$ by convention. Also, we let
$\fr(a):=0$ by convention for the empty $p$-tuple $a=()$ with
$p=0$~.
\begin{example}\label{35ex}
If $a=(3,5)$ then the gap set is $G(a)=\{1,2,4,7\}$~, and the
Frobenius number is $\fr(a)=7$~.
\end{example}
Classical results  of Schur and Sylvester, respectively,
assert that for all $p\geq 2$ and all $a=(a_1,\dots,a_p)$ with each
$a_i\geq 2$~, the Frobenius number obeys the upper bound
\begin{equation}\label{Frobenius-bound}
\fr(a)\, +\, 1 ~\leq~ \min\left\{(a_i-1)(a_j-1) ~:~ 1\leq i < j\leq
p\right\}\ ,
\end{equation}
with equality $\fr(a)+1=(a_1-1)(a_2-1)$ holding for $p=2$~. See
\cite{Bra} and references therein for proofs.

Define the \emph{restriction} of $M(a)$ by $\la\in \Z_+^p$ to be the
following subset of $M(a)$~:
$$M(a,\la) ~:=~ \{\mu a ~:~ \mu\in\Z_+^p\,,\ \mu\leq\la\}\ .$$

We start with a few simple facts.

\bp{Sym} For every $\la\in\Z_+^p$~, $M(a,\la)$ is symmetric on
$\{0,1,\dots,\la a\}$~, that is, we have that $g\in M(a,\la)$ if and
only if $\la a - g\in M(a,\la)$~. \ep

\bpr Indeed, $g=\mu a$ with $0\le \mu\le\la$ if and only if $\la a -
g=(\la-\mu) a$ with $0\le \la-\mu\le\la$~. \epr

Recall that for $z,s\in\Z$ and $Z\subseteq \Z$~, we let
$z+sZ:=\{z+sx ~:~ x\in Z\}$~.

\bp{GapObservation} For every $\la\in\Z_+^p$~, we have
\begin{equation}\label{GapEquation}
M(a,\la)\ \subseteq\ \{0,1,\dots,\la a\} \setminus \left(\, G(a)\cup
(\la a-G(a))\, \right)\ .
\end{equation}
\ep

\bpr Clearly, $M(a,\la)\ \subseteq\ \{0,1,\dots,\la a\} \setminus
G(a)$~. The claim now follows from Proposition \ref{Sym}. \epr

\bigskip
Call $\la\in\Z_+^p$ \emph{saturated for $a$} if (\ref{GapEquation})
holds for $\la$ with equality. In particular, if some $a_i=1$~, then
$\la$ saturated for $a$ implies $M(a,\la)=\{0,1,\dots,\la a\}$~.

\noindent {\bf Example \ref{35ex}, continued.} For $a=(3,5)$ and say
$\la=(3,4)$~, we have $\la a=29$, and it can be easily checked that
there are two values, namely $12=4\cdot 3 + 0\cdot 5$ and $17=4\cdot
3 + 1\cdot 5$~, that are not in $M(a,\la)$ but are in
$\{0,1,\dots,\la a\} \setminus \left(\, G(a)\cup (\la a-G(a))\,
\right)$~. Hence, in this case $\la$ is not saturated for $a$~.

Let $\max(a):=\max\{a_1,\dots,a_p\}$~.
Call $a=(a_1,\dots,a_p)$ \emph{divisible} if $a_i$ divides $a_{i+1}$
for $i=1,\dots p-1$~.
The following theorem asserts that, for any fixed primitive $a$~, every
(component-wise) sufficiently large $p$-tuple $\la$ is saturated for $a$~.

\bt{GapTheorem} Let $a=(a_1,\dots,a_p)$ be any primitive $p$-tuple.
Then the following statements hold:
\begin{enumerate}
\item Every $\la=(\la_1,\dots,\la_p)$ satisfying $\la_i\geq \max(a)$
for $i=1,\dots,p$ is saturated for $a$~.
\item For divisible $a$~, every $\la=(\la_1,\dots,\la_p)$ satisfying
$\la_i\geq{a_{i+1}\over{a_i}}-1$ for $i=1,\ldots,p-1$ is saturated for $a$~.
\end{enumerate}
\et

\bpr We begin with Part 1. As we go, we make
some claims  for which we employ somewhat tedious and lengthy
elementary arguments to carefully verify. We relegate proofs
of these claims, specifically Claim 1 and SubClaims 2.1--2.4, to the Appendix.

Suppose that $\lambda_i\ge\max(a)$~,
for $i=1,\ldots,p$~. Suppose that the result is false. Then there is
a $p$-tuple $\mu\in\Z^p_+$ so that $\mu a \leq  \lambda a$ but $\mu
a \notin M(a,\lambda)$~. By Proposition \ref{Sym}, we can assume
that $\mu a \le \frac{1}{2}\lambda a$~. Among all such $\mu$~,
choose one that has minimum \emph{violation}
$\sum_{i=1}^p(\mu_i-\la_i)^+$~.
Let $j$ be an index such that $\mu_j > \lambda_j$~.

\noindent \emph{Claim 1:} There are at least two indices $k$ for which
$\mu_k< \lambda_k/2$~.

Next, for every integer $ 0 \leq \gamma \leq a_j-1 $~,  consider the
two-variable  integer linear program:
\begin{alignat*}{1}
& \min\ x_l(\gamma) \nonumber \\
& \text{s.t. } a_j x_j(\gamma) - a_l x_l(\gamma)  =   \gamma a_k ~;\label{status}\tag*{$P_\gamma$}\\
& x_j(\gamma)~,~x_l(\gamma) \in \Z_+~.\nonumber
\end{alignat*}

\noindent \emph{Claim 2:} For some $\gamma \leq \lceil a_j/2\rceil$~,
there is a nonzero optimal solution to $P_\gamma$~,
such that $x_l(\gamma) \leq \lfloor a_j/2 \rfloor$~.

\noindent \emph{Proof of Claim 2:}
For the purpose of establishing Claim 2,
we assume, without loss of generality, that $\gcd(a_j,a_k,a_l)=1$~;
if this did not hold, we could just divide the integers
$a_j,a_k,a_l$ by their greatest common divisor, thus proving
a stronger result.

\medskip

\noindent \emph{SubClaim 2.1:} The  integer program \ref{status} is
feasible for all integers  $ 0 \leq \gamma ~(\leq a_j-1) $ that are
integer multiples of $\gcd(a_l,a_j)$~.

\noindent \emph{SubClaim 2.2:}
In fact,
for $\gamma = z_k \gcd(a_l,a_j)$ with $z_k \in \Z_+$~, we have that
$x^*_l(\gamma) =  z_l \gcd(a_k,a_j)$ for some $z_l \in \Z_+$~.

\noindent \emph{SubClaim 2.3:}
For $0\leq \gamma,\gamma'< a_j/\gcd(a_k,a_j)$~,
we have that $x^*_l(\gamma) \neq x^*_l(\gamma')$ for $\gamma \neq \gamma'$~.

\noindent \emph{SubClaim 2.4:} For integer $\gamma\ge a_j/ \gcd(a_k,a_j)$~,
we write $\gamma$ uniquely as
\[
\gamma ~=~ \gamma' ~+~ \mu a_j /  \gcd(a_k,a_j)~,
\]
with $\mu\in\Z_+$~, $\gamma'\in\Z_+$~, $\gamma'<a_j/\gcd(a_k,a_j)$~. Then we have that
\begin{eqnarray*}
x^*_l(\gamma') & = & x^*_l(\gamma)~,\\
x^*_j(\gamma') & = & x^*_j(\gamma) ~+~  \mu a_k /  \gcd(a_k,a_j)~.
\end{eqnarray*}

Now we are in position to complete the proof of Claim 2.
First, if $\gcd(a_l,a_j) \geq 2$~, then Claim 2 follows because
\begin{eqnarray*}
x_l(0) & := & a_j/\gcd(a_l,a_j)~,\\
x_j(0) & := & a_l/\gcd(a_l,a_j)
\end{eqnarray*}
is a feasible solution of $P_0$ with $x_l(0) \leq \lfloor a_j/2 \rfloor$~.
So, we can assume from now on that $\gcd(a_l,a_j)=1$~.

We denote by $\Omega$ the set of all integers $ 0 \leq \gamma \leq a_j-1$
for which \ref{status} is feasible.
 Next, assume that $\gcd(a_k,a_j)\geq 2$~.
 Then by what we have shown already,
 \[
  \left\{x^*_l(\gamma) ~:~ \gamma \in \Omega\right\} ~=~
  \left\{x^*_l(\gamma) ~:~ \gamma \in \Omega~,~ \gamma < a_j/\gcd(a_k,a_j)\right\}~.
  \]
Because $a_j/\gcd(a_k,a_j) \leq a_j/2$~, there is a
 $\gamma \leq a_j/\gcd(a_k,a_j) \leq a_j/2$ such
 that $P_\gamma$ has a feasible solution with
$x_l(\gamma)=1$~.
So we now can further assume that  $\gcd(a_k,a_j) = 1$~.

Then $x^*_l(\gamma) \neq x^*_l(\gamma')$ for all $\gamma \in \Omega$~,
$\gamma \neq \gamma'$
implies that
the cardinality of the set
     $\{x^*_l(\gamma) ~:~ 1 \leq \gamma \leq \lceil a_j/2\rceil\}$ is equal to
$\lceil a_j/2\rceil$~.
Because $x^*_l(\gamma)$ is an integer between 0 and $a_j - 1$~, it follows that there
must exist a $\gamma^*$  with $1 \leq \gamma^* \leq \lceil a_j/2\rceil$ such that
 $x^*_l(\gamma^*) \leq \lfloor a_j/2 \rfloor$~.
Hence we have established Claim 2.

Notice that this then also implies
that
\[
x_j^*(\gamma^*) ~a_j = \gamma^* a_k +
x_l^*(\gamma^*)~a_l \leq \max(a) \left( \gamma^* +  x_l^*(\gamma^*)
\right) \leq \max(a)~ a_j~,
\]
which implies $x_j^*(\gamma^*)
\leq \max(a)$~.

Now, define a new $p$-tuple $\nu$ by
\[
\nu_j:=\mu_j-x_j^*(\gamma^*)~,\quad
\nu_l:=\mu_l+x_l^*(\gamma^*)~,\quad \nu_k:=\mu_k+\gamma^*~, \quad
\mbox{and}\ \ \nu_i:=\mu_i\ \ \mbox{for all}\ \ i\neq j,k,l~.
\]

Because $x_j^*(\gamma^*)
\leq \max(a)$~, it follows that  $\nu_j> 0$~. Moreover, for
$i\in \{k,l\}$~, $0 \leq \nu_i \leq \la_i$~. Therefore $\nu$ is nonnegative,
satisfies $\nu a =\mu a= v$~, and has lesser violation than $\mu$~,
which is a contradiction to the choice of $\mu$~. So indeed $v\in M(a,\la)$~,
and we have established Part 1 of the theorem.

Before continuing, we note that a much simpler elementary argument
can be used to establish Part 1 of the theorem
under the stronger hypothesis: $\lambda_i\geq 2\max(a)$ for
$i=1,\ldots,p$~.

\vskip0.5cm
We next proceed with establishing Part 2 of the theorem. We
begin by using induction on $p$~. For $p=1$~, we have $a_1=1$~, and
every $\la=(\la_1)$ is saturated because every $0\leq v\leq \la
a=\la_1$ satisfies $v=\mu a=\mu_1$ for $\mu\leq\la$ given by
$\mu=(\mu_1)$ with $\mu_1=v$~.

Next consider $p>1$~. We use induction on $\la_p$~. Suppose first
that $\la_p=0$~. Let $a':=(a_1,\dots,a_{p-1})$ and
$\la':=(\la_1,\dots,\la_{p-1})$~. Consider any value $0\leq v\leq
\la a=\la'a'$~. Since $\la'$ is saturated by induction on $p$~,
there exists $\mu'\leq\la'$ with $v=\mu'a'$~. Then,
$\mu:=(\mu',0)\leq\la$ and $v=\mu a$~. So $\la$ is also saturated.
Next, consider $\la_p>0$~. Let
$\tau:=(\la_1,\dots,\la_{p-1},\la_p-1)$~. Consider any value $0\leq
v\leq \tau a=\la a-a_p$~. Since $\tau$ is saturated by induction on
$\la_p$~, there is a $\mu\leq\tau<\la$ with $v=\mu a$~, and so $v\in
M(a,\tau)\subseteq M(a,\la)$~. Moreover, $v+a_p={\hat \mu}a$ with
${\hat \mu}:=(\mu_1,\dots,\mu_{p-1},\mu_p+1)\leq\la$~, so $v+a_p\in
M(a,\la)$ as well. Therefore
\begin{equation}\label{divisible_equation}
\{0,1,\dots,\tau a\}\ \cup\ \{a_p,a_p+1,\dots,\la a\}\ \subseteq\
M(a,\la)\ .
\end{equation}
Now,
\[
\tau a\ =\ \sum_{i=1}^p\tau_i a_i\ \geq\ \sum_{i=1}^{p-1}\la_i a_i\
\geq \ \sum_{i=1}^{p-1}\left({{a_{i+1}}\over {a_i}}-1\right) a_i\ =
\ \sum_{i=1}^{p-1}(a_{i+1}-a_i)\ = \ a_p-1~,
\]
implying that the left-hand side of (\ref{divisible_equation}) is in
fact equal to $\{0,1,\dots,\la a\}$~. Therefore $\la$ is indeed
saturated. This completes the double induction, the proof of Part 2,
and the proof of the theorem.
\epr

\section{Obtaining an $r$-Best Solution}\label{Algorithm}

We can now combine all the ingredients developed in the previous
sections and provide our algorithm. Let $a=(a_1,\dots,a_p)$ be a
fixed primitive $p$-tuple. Define $\la=(\la_1,\dots,\la_p)$ by
$\la_i:=\max(a)$ for every $i$~. For $\mu\leq\la$ define
\[
I^\la_\mu\ :=\ \{i ~:~ \mu_i=\la_i\}\quad \hbox{  and  }\quad
a^\la_\mu\ :=\ \left({a_i\over \gcd(a_i ~:~ i\in I^\la_\mu)} ~:~
i\in I^\la_\mu\right)~.
\]
Finally, define
\begin{equation}\label{quality}
r(a)\quad :=\quad \sum_{\mu\leq\la} \fr(a^\la_\mu)~.
\end{equation}

\noindent The next corollary gives some estimates on $r(a)$~,
including a general bound implied by Theorem \ref{GapTheorem}.
\bc{Bound-Corollary} Let $a=(a_1,\dots,a_p)$ be any primitive
$p$-tuple. Then the following hold:
\begin{enumerate}
\item
 An upper bound on $r(a)$ is given by $r(a)\leq \left(2\max(a)\right)^p$~.
\item
For divisible $a$~, we have $r(a)=0$~.
\item
For $p=2$~, that is, for $a=(a_1,a_2)$~, we have $r(a)=\fr(a)$~.
\end{enumerate}
\ec
\bpr Define $\la=(\la_1,\dots,\la_p)$ by $\la_i:=\max(a)$ for
every $i$~. First note that if $I^\la_\mu$ is empty or a singleton
then $a^\la_\mu$ is empty or $a^\la_\mu=1$~, and hence
$\fr(a^\la_\mu)=0$~.

Part 1: As noted, $\fr(a^\la_\mu)=0$ for each $\mu\leq\la$ with
$|I^\la_\mu|\leq 1$~. There are at most  $2^p(\max(a))^{p-2}$
$p$-tuples $\mu\leq\la$ with $|I^\la_\mu|\geq 2$ and for each, the
bound of equation (\ref{Frobenius-bound}) implies
$\fr(a^\la_\mu)\leq (\max(a))^2$~. Hence
$$r(a)\ \leq \ 2^p(\max(a))^{p-2}(\max(a))^2\ \leq \
\left(2 \max(a)\right)^p\ .$$

Part 2: If $a$ is divisible, then the least entry of every nonempty
$a^\la_\mu$ is $1$~, and hence $\fr(a^\la_\mu)=0$ for every
$\mu\leq\la$~. Therefore $r(a)=0$~.

Part 3: As noted, $\fr(a^\la_\mu)=0$ for each $\mu\leq\la$ with
$|I^\la_\mu|\leq 1$~. For $p=2$~, the only $\mu\leq\la$ with
$|I^\la_\mu|=2$ is $\mu=\la$~. Because $a^\la_\la=a$~, we find that
$r(a)=\fr(a)$~. \epr

We are now in position to prove the following refined version of our
main theorem (Theorem \ref{Main}).
\vskip.2cm\noindent
\bt{MainR}
For every primitive $p$-tuple $a=(a_1,\dots,a_p)$~, with $r(a)$ as
in (\ref{quality}) above, there is an algorithm that, given any
independence system $S\subseteq\{0,1\}^n$ presented by a
linear-optimization oracle, weight vector
$w\in\{a_1,\dots,a_p\}^n$~, and function $f:\Z\rightarrow\R$
presented by a comparison oracle, provides an $r(a)$-best solution
to the nonlinear problem $\min\{f(wx) ~:~ x\in S\}$~, in time
polynomial in $n$~. Moreover:
\begin{enumerate}
\item
If $a_i$ divides $a_{i+1}$ for $i=1,\dots, p-1$~, then the algorithm
provides an optimal solution.
\item
For $p=2$~, that is, for $a=(a_1,a_2)$~, the algorithm provide an
$\fr(a)$-best solution.
\end{enumerate}
\et

\bpr Consider the following algorithm:

\begin{algorithm2e}[H]
\label{algo1}

{\bf input} independence system $S\subseteq\{0,1\}^n$ presented by a
linear-optimization oracle, $f:\Z\rightarrow\R$ presented by a
comparison oracle, and $w\in\{a_1,\ldots,a_p\}^n$~\;

{\bf define} $\la=(\la_1,\dots,\la_p)$ by $\la_i:=\max(a)$ for every
$i$~\;

\For {\emph{every choice of $p$-tuple} $\mu\in\Z^p_+~,~
\mu\leq\la$}{

{\bf find} by the algorithm of Lemma \ref{mu-lemma} an $x_\mu$
attaining $\max\{wx ~:~ x\in S^\la_\mu\}$ if any\;

{\bf if} $S^\la_\mu\neq\emptyset$ {\bf then} find by the algorithm
of Lemma \ref{strategy-lemma} an $x^*_\mu$ attaining $\min\{f(wx)
~:~ x\in\{0,1\}^n~,~ x\leq x_\mu\}$~\;

}

{\bf output} $x^*$ as one minimizing $f(wx)$ among the $x^*_\mu$~\ .

\end{algorithm2e}

\noindent First note that the number of $p$-tuples $\mu\leq\la$ and
hence of loop iterations and applications of the polynomial-time
algorithms of Lemma \ref{strategy-lemma} and Lemma \ref{mu-lemma} is
$\prod_{i=1}^p(\la_i+1)=(1+\max(a))^p$ which is constant since $a$
is fixed. Therefore the entire running time of the algorithm is
polynomial.

Consider any $p$-tuple $\mu\leq\la$ with $S^\la_\mu\neq\emptyset$~,
and let $x_\mu$ be an optimal solution of $\max\{wx ~:~ x\in
S^\la_\mu\}$ determined by the algorithm. Let $I:=I^\la_\mu=\{i ~:~
\mu_i=\la_i\}$~, let $g:=\gcd(a_i ~:~ i\in I)$~, let ${\bar
a}:=a^\la_\mu={1\over g}(a_i ~:~ i\in I)$~, and let $h:=\sum\{\mu_i
a_i ~:~ {i\notin I}\}$~. For each point $x\in\{0,1\}^n$ and for each
$i=1,\dots,p$~, let as usual $\la_i(x):=|\supp(x)\cap N_i|$~, where
$N_i=\{j ~:~ w_j=a_i\}$~, and let ${\bar\la}(x):=(\la_i(x) ~:~ i\in
I)$~. By the definition of $S^\la_\mu$ in equation
(\ref{partition-equation}) and of $I$ above, for each $x\in
S^\la_\mu$ we have
$$wx\ =\ \sum_{i\notin I}\la_i(x) a_i + \sum_{i\in I}\la_i(x) a_i
\ =\ \sum_{i\notin I}\mu_i a_i + g \sum_{i\in I}\la_i(x) {1\over
g}a_i \ =\ h + g {\bar \la}(x){\bar a}\ .$$ In particular, for every
$x\in S^\la_\mu$ we have $wx\in h+g M({\bar a})$ and $wx\leq
wx_\mu=h+g{\bar\la}(x_\mu){\bar a}$~, and therefore
$$w\cdot S^\la_\mu\ \subseteq
\ h\ +\ g\left(M({\bar a})\cap\{0,1\dots,{\bar\la}(x_\mu){\bar
a}\}\right)\ .$$ Let $T:=\{x ~:~ x\leq x_\mu\}$~. Clearly, for any
${\bar\nu}\leq{\bar\la}(x_\mu)$ there is an $x\in T$ obtained by
zeroing out suitable entries of $x_\mu$ such that
${\bar\la}(x)={\bar\nu}$ and $\la_i(x)=\la_i(x_\mu)=\mu_i$ for
$i\notin I$~, and hence $wx=h + g{\bar\nu}{\bar a}$~. Therefore
$$h\ +\ g M\left({\bar a},{\bar\la}(x_\mu)\right)\ \subseteq\ w\cdot T\ .$$
Since $x_\mu\in S^\la_\mu$~, by the definition of $S^\la_\mu$ and
$I$~, for each $i\in I$ we have
$$\la_i(x_\mu)\ =\ |\supp(x)\cap N_i|\ \geq\ \mu_i\ =\ \la_i
\ =\ \max(a)\ \geq \max({\bar a})\ .$$ Therefore, by Theorem
\ref{GapTheorem}, we conclude that ${\bar\la}(x_\mu)=(\la_i(x_\mu)
~:~ i\in I)$ is saturated for $\bar a$ and hence
$$M\left({\bar a},{\bar\la}(x_\mu)\right)\ =\
\left(M({\bar a})\cap\{0,1\dots,{\bar\la}(x_\mu){\bar a}\}\right)\
\setminus \left({\bar\la}(x_\mu){\bar a}-G({\bar a})\right)\ .$$
This implies that
\[
w\cdot S^\la_\mu\ \setminus\ w\cdot T\ \subseteq\ h\ +\ g
\left({\bar\la}(x_\mu){\bar a}-G({\bar a})\right)~,
\]
and hence
\[
|w\cdot S^\la_\mu\ \setminus\ w\cdot T|\ \leq\ |G({\bar a})|\ =\
\fr({\bar a})~.
\]
Therefore, as compared to the objective value of the optimal
solution $x^*_\mu$ of
\[
\min\{f(wx) ~:~ x\in T\}=\min\{f(wx) ~:~ x\leq x_\mu\}
\]
determined by the algorithm, at  most $\fr({\bar a})$ better
objective values are attained by points in  $S^\la_\mu$~.

Since $S=\biguplus_{\mu\leq\la} S^\la_\mu$ by Proposition
\ref{Partition}, the independence system $S$ has altogether at most
\[
\sum_{\mu\leq\la}\fr(a^\la_\mu)\ =\ r(a)
\]
better objective values $f(wx)$ attainable than that of the solution
$x^*$ output by the algorithm. Therefore $x^*$ is indeed an
$r(a)$-best solution to the nonlinear optimization problem over the
(singly) weighted independence system. \epr

In fact, as the above proof of Theorem \ref{MainR} shows, our algorithm provides
a better, $g(a)$-best, solution, where $g(a)$ is defined as follows in terms of
the cardinalities of the gap sets of the subtuples $a^\la_\mu$ with $\lambda$ defined
again by $\lambda_i:=2\max(a)$ for all $i$ (in particular,  $g(a)=|G(a)|$ for $p=2$),
\begin{equation}\label{tighter_quality}
g(a)\quad :=\quad \sum_{\mu\leq\la} |G(a^\la_\mu)|~.
\end{equation}

\section{Finding an Optimal Solution Requires Exponential Time}\label{LowerBound}

We now demonstrate that our results are best possible in the
following sense. Consider $a:=(2,3)$. Because $F(2,3)=1$, Theorem \ref{Main} (Part 2)
assures that our algorithm produces a $1$-best solution in
polynomial time. We next establish a refined version of Theorem \ref{Main2},
showing that a $0$-best (i.e., optimal)
solution {\em cannot} be found in polynomial time.

\bt{lower_bound} There is no polynomial time algorithm for computing
a $0$-best (i.e., optimal) solution of the nonlinear optimization
problem $\min\{f(wx)\,:\,x\in S\}$ over an independence system
presented by a linear optimization oracle with $f$ presented by a
comparison oracle and weight vector $w\in\{2,3\}^n$. \break In fact,
to solve the nonlinear optimization problem over every independence
system $S$ with a ground set of $n=4m$ elements with $m\geq 2$, at
least ${2m\choose m+1}\geq 2^m$ queries of the oracle presenting $S$
are needed.
\et

\bpr Let $n:=4m$ with $m\geq 2$,
$I:=\{1,\dots,2m\}$, $J:=\{2m+1,\dots,4m\}$, and let $w:=2\cdot{\bf
1}_I+3\cdot{\bf 1}_J$~. For  $E\subseteq \{1,\dots,n\}$ and any
nonnegative integer $k$~, let $E\choose k$ be the set of all $k$-element subsets of $E$.
For $i=0,1,2$~, let
\[
T_i ~:=~ \left\{x={\bf 1}_A+{\bf 1}_B  ~:~
A\in {I\choose {m+i}}\,,\ B\in {J\choose {m-i}}\right\}\ \subset\
\{0,1\}^n~.
\]
 Let $S$ be the independence system generated by
$T_0\cup T_2$, that is,
\[
S ~:=~
\left\{z\in\{0,1\}^n ~:~ z\leq x~,~ \mbox{for some}~ x\in T_0\cup T_2\right\}~.
\]
Note that the $w$-image of $S$ is
\[
w\cdot S=\{0,\dots,5m\}\setminus\{1,5m-1\}~.
\]

For every $y\in T_1$~, let $S_y:=S\cup\{y\}$~. Note that each $S_y$ is
an independence system as well, but with $w$-image
\[
w\cdot S_y=\{0,\dots,5m\}\setminus\{1\}~;
\]
that is, the $w$-image of each $S_y$ is precisely the $w$-image of $S$ augmented
by the value $5m-1$~.

Finally, for each vector $c\in\Z^n$~, let
\[
Y(c) ~:=~ \left\{y\in T_1\,:\, cy>\max\{cx:x\in S\}\right\}~.
\]

\noindent \emph{Claim:}
$|Y(c)|\leq{2m\choose m-1}$ for every $c\in\Z^n$~.

\noindent \emph{Proof of Claim:}
Consider two elements (if any) $y,z\in Y(c)$~. Then
$y={\bf 1}_A+{\bf 1}_B$ and $z={\bf 1}_U+{\bf 1}_V$ for some
$A,U\in{I\choose m+1}$ and $B,V\in{J\choose m-1}$. Suppose,
indirectly, that $A\neq U$ and $B\neq V$. Pick $a\in A\setminus U$
and $v\in V\setminus B$~.
Consider the following vectors,
\begin{eqnarray*}
x^0 &:=& y-{\bf 1}_a+{\bf 1}_v\ \in\ T_0 ~, \\
x^2 &:=& z+{\bf 1}_a-{\bf 1}_v\ \in\ T_2 ~.
\end{eqnarray*}
 Now $y,z\in Y(c)$ and
$x^0,x^2\in S$ imply the contradiction
\begin{eqnarray*}
c_a-c_v=cy-cx^0>0~, \\
c_v-c_a=cz-cx^2>0~.
\end{eqnarray*}
This implies that all vectors in $Y(c)$ are of the form ${\bf
1}_A+{\bf 1}_B$ with either $A\in{I\choose m+1}$ fixed, in which
case $|Y(c)|\leq{2m\choose m-1}$, or $B\in{J\choose m-1}$ fixed, in
which case $|Y(c)|\leq{2m\choose m+1}={2m\choose m-1}$~, as claimed.

Continuing with the proof of our theorem, consider any algorithm, and let
$c^1,\dots,c^p\in\Z^n$ be the sequence of oracle queries made by the
algorithm. Suppose that $p<{2m\choose m+1}$~. Then
\[
\left|\bigcup_{i=1}^p Y(c^i)\right|\ \leq\ \sum_{i=1}^p |Y(c^i)|
\ \leq\ p{2m\choose m-1}\ <\ {2m\choose m+1}{2m\choose m-1}\ =\
|T_1|~.
\]
 This implies that there exists some $y\in T_1$ that
is an element of none of the $Y(c^i)$~, that is, satisfies $c^iy\leq
\max\{c^ix:x\in S\}$ for each $i=1,\dots,p$~. Therefore, whether the
linear optimization oracle presents $S$ or $S_y$~, on each query
$c^i$ it can reply with some $x^i\in S$ attaining
\[
c^ix^i\ =\ \max\{c^ix:x\in S\}\ =\ \max\{c^ix:x\in S_y\}~.
\]
Therefore, the algorithm cannot tell whether the oracle presents $S$
or $S_y$ and hence can neither compute the $w$-image of the independence
system nor solve the nonlinear optimization problem
correctly.
\epr

\section{Discussion}\label{Conclusion}

We view this article as a first step in understanding the complexity
of the general nonlinear optimization problem over an independence
system presented by an oracle. Our work raises many intriguing
questions including the following. Can the saturated $\la$ for $a$
be better understood or even characterized? Can a saturated $\la$
smaller than that with $\la_i=\max(a)$ be determined for every $a$
and be used to obtain better running-time guarantee for the
algorithm of Theorem \ref{Main} and better approximation quality
$r(a)$~? Can tighter bounds on $r(a)$ in equation (\ref{quality})
and $g(a)$ in equation (\ref{tighter_quality}) and possibly formulas
for $r(a)$ and $g(a)$ for small values of $p$, in particular $p=3$, be derived?
For which primitive $p$-tuples $a$ can an exact solution to the nonlinear
optimization problem over a (singly) weighted independence system be
obtained in polynomial time, at least for small $p$, in particular
$p=2$~? For $p=2$ we know that we can when $a_1$ divides $a_2$~, and
we cannot when $a:=(2,3)$~, but we do not have a complete characterization.
How about $d=2$~? While this includes the notorious exact
matching problem as a special case, it may still be that a
polynomial-time solution is possible. And how about larger, but
fixed, $d$~?

In another direction, it can be interesting to consider the problem for functions
$f$ with some structure that helps to localize minima. For instance,
if $f:\R\rightarrow\R$ is concave or even more generally
quasiconcave (that is, its ``upper level sets'' $\{z\in\R ~:~ f(z)
\ge \tilde{f}\}$ are convex subsets of $\R$~, for all $\tilde{f}\in
\R$~; see \cite{ADSZ}, for example), then the optimal value
$\min\{f(wx)\, :\, x \in S\}$ is always attained on the boundary of
${\rm conv}(w\cdot S)$~, i.e., if $x^*$ is a minimizer, then either
$wx^*=0$ or $wx^*$ attains $\max\{wx \, : \, x\in S\}$~, so the
problem is easily solvable by a single query to the
linear-optimization oracle presenting $S$~ and a single query to the
comparison oracle of $f$~. Also, if $f$ is convex or even more
generally quasiconvex (that is, its ``lower level sets'' $\{z\in\R
~:~ f(z) \le \tilde{f}\}$ are convex subsets of $\R$~, for all
$\tilde{f}\in \R$), then a much simplified version of the algorithm
(from the proof of Theorem \ref{MainR}) gives an $r$-best solution
as well, as follows.

\bp{Convex} For every primitive $p$-tuple $a=(a_1,\dots,a_p)$~,
there is an algorithm that, given independence system
$S\subseteq\{0,1\}^n$ presented by a linear-optimization oracle,
weight vector $w\in\{a_1,\dots,a_p\}^n$~, and quasiconvex function
$f:\R\rightarrow\R$ presented by a comparison oracle, provides a
$(\max(a)-1)$-best solution to the nonlinear problem
$\min\{f(wx)~:~x\in S\}$~, in time polynomial in $n$~. \ep

\bpr We could describe the construction as a specialization of the
algorithm from the proof of Theorem \ref{MainR}, but it is more
clear to just present it directly. We first use our
linear-optimization oracle to find $x^*$ attaining $\max\{wx \, : \,
x\in S\}$~. Then, by repeatedly, and in an arbitrary order,
decreasing a single component of the point by unity, we obtain a
sequence of points
\[
x^k:=x^*~\ge~x^{k-1}~\ge~\ldots~\ge~x^0:={\bf 0}~,
\]
with $k=\sum_{j=1}^n x_j^*\le n$~. Let $\breve{f}:=\min\{f(wx^t) ~:~
0\le t\le k\}$~.

Next, using the comparison oracle (a linear number of times), we
find the least and greatest indices $t$~, say $t_{\min}$ and
$t_{\max}$ respectively, for which $x^t$ minimizes $f(wx^t)$~.
Quasiconvexity of $f$ implies that
\[
f(wx^{t})~=~\breve{f}~, \mbox{ for } t_{\min}\le t\le t_{\max}~.
\]
Moreover, quasiconvexity implies that there is an index $s$~,
satisfying $t_{\min}-1 \le s \le t_{\max}$~, such that \emph{all}
points $z\in[0,wx^*]\cap \Z$ having $f(z)<\breve{f}$ are in
$[wx^s+1,wx^{s+1}-1]\cap \Z$ (that is, in \emph{one} of the
$t_{\max}-t_{\min}+2$ intervals $[wx^{t},wx^{t+1}]$ beginning with
the one immediately to the left of $t_{\min}$ and ending with the
one immediately to the right of $t_{\max}$ --- and not the endpoints
of that interval).

The result now follows by noticing that
\[
wx^{t+1} - wx^{t} \leq \max(a)~,\quad \mbox{for}~ t=0,\ldots,k-1~,
\]
in particular for $t=s$~.
\epr

In yet another direction, it would be interesting to consider other (weaker or
stronger) oracle presentations of the independence system $S$.
While a membership oracle suffices for nonlinear
optimization when $S$ is a matroid \cite{BLMORWW}, in general
it is much too weak, as the following proposition shows.

\bp{membership}
There is no polynomial time algorithm for solving the nonlinear
optimization problem $\min\{f(wx)\,:\,x\in S\}$ over an independence
system presented by a membership oracle with $f$ presented by a
comparison oracle, even with all weights equal to $1$, that is,
for $p=1$, $a=1$, $w=(1,\dots,1)$.
\ep

\bpr
Let $n:=2m$~, let $w:=\sum_{i=1}^n{\bf 1}_i=(1,\dots,1)$~, and let
\[
S\ :=\ \{x\in\{0,1\}^n\ :\ \supp(x)\leq m-1\}~.
\]
For each $y\in\{0,1\}^n$ with $\supp(y)=m$~, let $S_y:=S\cup\{y\}$~.
Note that
\[
w\cdot S=\{0,1,\dots,m-1\}\,,\quad\quad
w\cdot S_y=\{0,1,\dots,m-1,m\}~.
\]
Now, suppose an algorithm queries the membership oracle less than
$n\choose m$ times. Then some $y\in\{0,1\}^n$ with $\supp(y)=m$ is not
queried, and so the algorithm cannot tell whether the oracle presents
$S$ or $S_y$ and hence can neither compute the image nor solve the
nonlinear optimization problem correctly.
\epr

\section*{Acknowledgment}

This research was supported by the Mathematisches Forschungsinstitut
Oberwolfach during a stay within the Research in Pairs Programme.

\section*{Appendix}

\noindent \emph{Claim 1:} There are at least two indices $k$ for which
$\mu_k< \lambda_k/2$~.

\noindent \emph{Proof of Claim 1:}
We note that $\mu_j < \lambda_j$ trivially implies
\begin{equation}\label{jeq1}
0 \le a_j \left(\mu_j -\lambda_j-1\right)~.
\end{equation}
Also, $\mu a \le \frac{1}{2}\lambda a$ can be written as
\begin{equation}\label{jeq2}
\sum_{k\not= j} a_k \left(\mu_k - \lambda_k/2\right) \le
a_j \left(-\mu_j+ \lambda_j/2\right)~.
\end{equation}
Now, adding (\ref{jeq1}) and (\ref{jeq2}), we obtain
\begin{equation}\label{jeq3}
\sum_{k\not= j} a_k \left(\mu_k - \lambda_k/2\right) \le
-a_j \left(\lambda_j/2 + 1\right)~.
\end{equation}
The right-hand side of (\ref{jeq3}) is negative, therefore the
left-hand side must also be negative.
Suppose that there is but a single index $k$ for which a summand
on the left-hand side of (\ref{jeq3}) is negative. Then, we have
\[
a_k \left(\mu_k - \lambda_k/2\right) \le
-a_j \left(\lambda_j/2 + 1\right)~,
\]
which implies
\begin{equation}\label{jeq4}
\max(a) \left(\mu_k - \lambda_k/2\right) \le
-a_j \left(\max(a)/2 + 1\right)~.
\end{equation}

We observe that we must have $\mu_k - \lambda_k >
-a_j$~, otherwise we could decrease the violation by decreasing
$\mu_j$ by $a_k$ and increasing $\mu_k$ by $a_j$~.
But $\mu_k - \lambda_k > -a_j$ implies that
\begin{equation}\label{jeq5}
\mu_k - \lambda_k/2 \ge
-a_j + 1 + \lambda_k/2 \ge
-a_j + 1 + \max(a)/2
~.
\end{equation}

Next, we combine (\ref{jeq4}) and  (\ref{jeq5}) to arrive at
\[
\max(a)\left( -a_j + 1 + \max(a)/2 \right) \le
-a_j \left(\max(a)/2 + 1\right)~,
\]
or, equivalently,
\[
a_j \left( \max(a)/2-1~\right) \ge \max(a) \left(\max(a)/2+1\right),
\]
which cannot hold.

So Claim 1 is established.

\medskip

\noindent \emph{SubClaim 2.1:} The  integer program \ref{status} is
feasible for all integers  $ 0 \leq \gamma ~(\leq a_j-1) $ that are
integer multiples of $\gcd(a_l,a_j)$~.

\noindent \emph{Proof of SubClaim 2.1:}
Suppose that $\gamma:=z_k~ \gcd(a_l,a_j)$~, for some $z_k\in\Z_+$~.

By B\'ezout's Lemma, there are integers $\beta_j~,~\beta_l$
such that
\[
a_j \beta_j + a_l \beta_l = \gcd(a_l,a_j)~.
\]
Moreover, there is an infinite family indicated by
\[
a_j \biggl(\beta_j + t a_l/\gcd(a_l,a_j)\biggr) + a_l \biggl( \beta_l - t a_j/\gcd(a_l,a_j)\biggr) = \gcd(a_l,a_j)~,
\]
with $t$ ranging over $\Z$~.

Multiplying through by $z_k a_k$~, and rearranging terms, we obtain
\begin{eqnarray*}
a_j  \left(z_k a_k \biggl(\beta_j + t a_l/\gcd(a_l,a_j)\biggr) \right)
+ a_l \left(z_k a_k \biggl( \beta_l - t a_j/\gcd(a_l,a_j)\biggr)\right)   &= &z_k \gcd(a_l,a_j) a_k\\
&=& \gamma a_k~.
\end{eqnarray*}

Now, for a sufficiently large positive integer $t$~, we will have
\[
\beta_l - t a_j/\gcd(a_l,a_j) ~\le~ 0~,
\]
and so
\begin{eqnarray*}
x_j(\gamma) &:=& z_k a_k \biggl(\beta_j + t a_l/\gcd(a_l,a_j)\biggr)~;\\
-x_l(\gamma) &:=& z_k a_k \biggl( \beta_l - t a_j/\gcd(a_l,a_j)\biggr)
\end{eqnarray*}
will be a feasible solution to \ref{status}~.
Thus we have established SubClaim 2.1.

\medskip

\noindent \emph{SubClaim 2.2:}
In fact,
for $\gamma = z_k \gcd(a_l,a_j)$ with $z_k \in \Z_+$~, we have that
$x^*_l(\gamma) =  z_l \gcd(a_k,a_j)$ for some $z_l \in \Z_+$~.

\noindent \emph{Proof of SubClaim 2.2:}
\begin{eqnarray*}
a_l x^*_l(\gamma) & = & a_jx^*_j(\gamma) ~-~ \gamma a_k\\
                  & = & \biggl(a_jx^*_j(\gamma)/\gcd(a_k,a_j) ~-~ \gamma a_k/\gcd(a_k,a_j)  \biggr) \gcd(a_k,a_j)~.
\end{eqnarray*}
As $\gcd(a_k,a_j)$ divides both $a_j$ and $a_k$~, we have
\[
a_lx^*_l(\gamma) = z~ \gcd(a_k,a_j)~,
\]
for some $z\in\Z_+$~, and hence
\[
x^*_l(\gamma) = (z/a_l) \gcd(a_k,a_j)~.
\]
As $\gcd(a_l,\gcd(a_k,a_j))=1$~, it is clear that $a_l$ must divide $z$ (after all
$x^*_l(\gamma)\in\Z)$, and hence SubClaim 2.2 is established.

\medskip

\noindent \emph{SubClaim 2.3:}
For $0\leq \gamma,\gamma'< a_j/\gcd(a_k,a_j)$~,
we have that $x^*_l(\gamma) \neq x^*_l(\gamma')$ for $\gamma \neq \gamma'$~.

\noindent \emph{Proof of SubClaim 2.3:}
Suppose the contrary.
 Without loss of generality,
 $\gamma' > \gamma$~.
 Then we have the following two equations:
\begin{eqnarray}
a_k\gamma' + a_l x^*_l(\gamma') &=& x^*_j(\gamma')a_j~;\label{one}\\
 a_k \gamma + a_l x^*_l(\gamma ) &=& x^*_j(\gamma)a_j~.\label{two}
\end{eqnarray}
Subtracting (\ref{two}) from (\ref{one}) gives
\begin{eqnarray}
  a_k (\gamma'-\gamma)&=&  a_j\left(x^*_j(\gamma')-x^*_j(\gamma)\right)~.\label{three}
\end{eqnarray}
Because $\gamma'>\gamma$~, the left-hand side  of (\ref{three}) is positive,
which implies that  $x^*_j(\gamma')-x^*_j(\gamma) > 0$~.

But $\gamma'-\gamma < a_j/\gcd(a_k,a_j)$~.
 This contradicts that $\gcd(a_j/\gcd(a_k,a_j) , a_k/\gcd(a_k,a_j)) = 1$~,
because every positive integer solution of
$a_k x_k = a_j x_j$
is a  positive multiple of
\begin{eqnarray*}
x_k&:=&a_j/\gcd(a_k,a_j)~,\\
x_j&:=&a_k/\gcd(a_k,a_j))~.
\end{eqnarray*}
Thus we have established SubClaim 2.3.

\medskip

\noindent \emph{SubClaim 2.4:} For integer $\gamma\ge a_j/ \gcd(a_k,a_j)$~,
we write $\gamma$ uniquely as
\[
\gamma ~=~ \gamma' ~+~ \mu a_j /  \gcd(a_k,a_j)~,
\]
with $\mu\in\Z_+$~, $\gamma'\in\Z_+$~, $\gamma'<a_j/\gcd(a_k,a_j)$~. Then we have that
\begin{eqnarray*}
x^*_l(\gamma') & = & x^*_l(\gamma)~,\\
x^*_j(\gamma') & = & x^*_j(\gamma) ~+~  \mu a_k /  \gcd(a_k,a_j)~.
\end{eqnarray*}

\noindent \emph{Proof of SubClaim 2.4:}
We can directly check feasibility:
\[
a_j \left( x^*_j(\gamma) ~+~  \mu a_k /  \gcd(a_k,a_j) \right) ~+~ a_l x^*_l(\gamma)
   ~=~ \left( \gamma' ~+~ \mu a_j /  \gcd(a_k,a_j) \right) a_k~.
\]
Moreover, there is no feasible
solution $\bar{x}$ for  $P_{\gamma'}$ having
$\bar{x}_l < x^*_l(\gamma)$~, because if there were, we would simply add
$\mu a_k /  \gcd(a_k,a_j)$ to $\bar{x}_j$~, and leave $\bar{x}_l$ unchanged,
to produce a feasible solution
for $P_{\gamma}$ having objective value less than $x^*_l(\gamma)$~, a
contradiction. Thus we have established SubClaim 2.4.

\vskip.6cm\noindent {\small Jon Lee}\newline \emph{IBM T.J. Watson
Research Center, Yorktown Heights, NY 10598, USA}\newline
\emph{email: jonlee{\small @}us.ibm.com}, \ \
\emph{http://www.research.ibm.com/people/j/jonlee}

\vskip.3cm\noindent {\small Shmuel Onn}\newline \emph{Technion -
Israel Institute of Technology, 32000 Haifa, Israel}\newline
\emph{email: onn{\small @}ie.technion.ac.il}, \ \
\emph{http://ie.technion.ac.il/{\small $\sim$onn}}

\vskip.3cm\noindent {\small Robert Weismantel}\newline
\emph{Otto-von-Guericke Universit\"at Magdeburg, D-39106 Magdeburg,
Germany}\newline \emph{email: weismantel{\small
@}imo.math.uni-magdeburg.de}, \ \
\emph{http://www.math.uni-magdeburg.de/{\small $\sim$weismant}}

\end{document}